\documentclass[a4paper,12pt]{article}
\usepackage{amssymb}
\usepackage{amsmath}
\usepackage{amsthm}
\usepackage{graphicx}
\numberwithin{equation}{section}
\newcommand{\ep}{\varepsilon}
\newcommand{\la}{\lambda}
\newcommand{\va}{\varphi}
\newcommand{\ppp}{\partial}
\newcommand{\ddd}{\mbox{div}\thinspace}

\newcommand{\weight}{e^{2s\va}}

\newcommand{\R}{\mathbb{R}}

\newcommand{\N}{\mathbb{N}}

\newcommand{\www}{\widetilde}

\newcommand{\ssss}{\sum_{i,j=1}^n}
\newcommand{\oooo}{\overline}

\newcommand{\XXX}{(x_1,x_2,\gamma(x_1,x_2))}

\allowdisplaybreaks

\title{Carleman estimate for the Navier-Stokes equations and
an application to a lateral Cauchy problem}

\author{Mourad Bellassoued\\
Department of Mathematics, Faculty of Sciences of Bizerte\\
University of  Carthage, 7021 Jarzouna Bizerte, Tunisia\\
e-mail: {\tt mourad.bellassoued@fsb.rnu.tn}\\
Oleg Imanuvilov\\
Department of Mathematics, Colorado State University\\
101 Weber Building, Fort Collins, CO 80523-1874, USA\\
e-mail: {\tt oleg@math.colostate.edu}\thanks{Partially supported
by NSF grant DMS 1312900}\\
Masahiro Yamamoto\\
Graduate School of Mathematical Sciences, the University of Tokyo\\
3-8-1 Komaba, Meguro-ku, Tokyo, 153-8914, Japan.\\
e-mail: {\tt myama@ms.u-tokyo.ac.jp}\thanks
{Partially supported by 
Grant-in-Aid for Scientific Research (S) 15H05740 
of Japan Society for the Promotion of Science}
}

\date{}

\begin{document}
 \maketitle


\begin{abstract}
We consider the nonstationary linearized
Navier-Stokes equations in a bounded domain and
first we prove a Carleman estimate with a regular weight function.
Second we apply the Carleman estimate to a lateral Cauchy problem
for the Navier-Stokes equations and prove the H\"older stability
in determining the velocity and pressure field in an interior
domain.
\end{abstract}

\baselineskip 18pt
\section{Introduction}

Let $\Omega \subset \R^n$, $n=2,3$, be a bounded domain
with smooth boundary $\ppp\Omega$ (e.g., of $C^2$-class),
and let $\nu=\nu(x)$ be the outward unit
normal vector on $\ppp\Omega$ at $x$.
We set $Q := \Omega \times (0,T)$.

We consider the linearized Navier-Stokes equations for
an incompressible viscous fluid:
$$
\ppp_t v(x, t) - \kappa\Delta v(x, t) + (A\cdot\nabla)v
+ (v\cdot \nabla)B + \nabla p = F(x, t) \quad
\mbox{in $Q$},                                           \eqno{(1.1)}
$$
and
$$
\ddd v(x, t) = 0 \quad \mbox{in $Q$}.      \eqno{(1.2)}
$$
Here $v = (v_1, \cdots, v_n)^T$, $n=2,3$, $\cdot^T$ denotes the transpose
of matrices, $\kappa > 0$ is a constant describing the viscosity,
and for simplicity we assume that the density is one.
Let $\ppp_t = \frac{\ppp}{\ppp t}$, $\ppp_j = \frac{\ppp}{\ppp x_j}$,
$1\le j\le n$, $\Delta = \sum_{j=1}^n \ppp_j^2$, $\nabla =
(\ppp_1,\cdots,\ppp_n)^T$, $\nabla_{x,t} = (\nabla, \ppp_t)^T$,
$\ppp_x^{\beta} = \ppp_1^{\beta_1}\cdots \ppp_n^{\beta_n}$ with
$\beta= (\beta_1, \cdots, \beta_n) \in (\N \cup \{0\})^n$,
$\vert \beta\vert = \beta_1 + \cdots + \beta_n$,
$$
(w\cdot \nabla)v = \left( \sum^n_{j=1} w_j\ppp_jv_1, \cdots,
\sum^n_{j=1} w_j\ppp_jv_n\right)
^T
$$
for $v = (v_1,\cdots,v_n)^T$ and $w = (w_1,\cdots,w_n)^T$.
Throughout this paper, we assume
$$
A \in W^{2,\infty}(Q), \quad
\nabla B \in L^{\infty}(Q).                          \eqno{(1.3)}
$$

In this paper, we establish a Carleman estimate with a regular
weight function and apply it to a lateral Cauchy problem for the 
Navier-Stokes equations and prove the H\"older stability in 
an arbitrarily given interior domain.
For stating the main results, we introduce notations.
Let $I_n$ be the $n\times n$ identity matrix and let the
stress tensor $\sigma(v,p)$ be defined by the $n\times n$ matrix
$$
\sigma(v,p) := \kappa (\nabla v + (\nabla v)^T) - pI_n,
$$ where $\kappa$ is some positive constant.
We assume
$$
d \in C^2(\overline{\Omega}), \quad \vert \nabla d(x)\vert > 0\quad
\mbox{on $\overline{\Omega}$}                  \eqno{(1.4)}
$$
and we arbitrarily choose $t_0 \in (0,T)$ and $\beta > 0$.  We set
$$
\psi(x,t) = d(x) - \beta(t-t_0)^2, \quad
\va(x,t) = e^{\lambda\psi(x,t)}
$$
with a sufficiently fixed large constant $\lambda>0$.
We choose a non-empty relatively open subboundary $\Gamma \subset \ppp\Omega$
arbitrarily.

Let $D \subset Q$ be a bounded domain with smooth boundary $\ppp Q$ such that
$\overline{D\cap (\ppp\Omega \times (0,T))} \subset \Gamma \times (0,T)$.

For $k, \ell \in \N \cup \{0\}$, we set
$$
H^{k,\ell}(D) = \{ v \in L^2(D); \thinspace
\ppp_x^{\beta}v \in L^2(D), \thinspace \vert \beta\vert \le k,
\thinspace \ppp_t^{j}v \in L^2(D) \quad 0\le j\le \ell\}
$$
and
$$
\Vert (v,p)\Vert^2_{ \mathcal X_s(D)}:=
\int_D \left\{ \frac{1}{s^2}\left( \vert \ppp_tv\vert^2
+ \sum_{i,j=1}^n \vert \ppp_i\ppp_jv\vert^2\right)
+ \vert \nabla v\vert^2 + s^2\vert v\vert^2
+ \frac{1}{s}\vert \nabla p\vert^2 + s\vert p\vert^2
\right\} \weight dxdt.
$$

We are ready to state our Carleman estimate.
\\
{\bf Theorem 1}.\\
{\it There exist constants $s_0 > 0$ and $C>0$, independent of $s$, such that
$$
\Vert (v,p)\Vert^2_{\mathcal X_s(D)}
\le C\int_D \vert F\vert^2 \weight dxdt
+ C\int_D (\vert h\vert^2 + \vert \nabla_{x,t} h\vert^2)\weight dxdt
$$
$$
+ Ce^{Cs}(\Vert v\Vert^2_{L^2(0,T;H^{\frac{3}{2}}(\Gamma))}
+ \Vert \ppp_tv\Vert^2_{L^2(0,T;H^{\frac{1}{2}}(\Gamma))}
+ \Vert \sigma(v,p)\nu\Vert^2_{L^2(0,T;H^{\frac{1}{2}}(\Gamma))})
                                              \eqno{(1.5)}
$$
for all $s \ge s_0$ and $(v,p) \in H^{2,1}(D) \times H^{1,0}(D)$
satisfying (1.1),
$$
\ddd v = h \quad \mbox{in $D$ with $h \in H^{1,1}(D)$},
$$
and
$$
\left\{
\begin{array}{rl}
& v(\cdot,0) = v(\cdot,T) = 0
\quad \mbox{in $\Omega$},\\
& \vert v\vert = \vert \nabla v\vert = \vert p\vert
= 0 \quad \mbox{on $\ppp D \setminus
(\Gamma \times (0,T))$}.
\end{array}
\right.                                   \eqno{(1.6)}
$$}

This is a Carleman estimate for the linearized Navier-Stokes
equations (1.1) with (1.2) with boundary data on $\Gamma \subset
\ppp\Omega$.  

Boulakia \cite{Bo} proves a Carleman estimate with a weight function
similar to ours for the homogeneous Stokes equations:
$\ppp_tv = \Delta v - \nabla p$ and $\ddd v = 0$ with extra
interior or boundary data.
The Carleman estimate in \cite{Bo} requires a stronger norm of boundary
data than our Carleman estimate if it is applied to
the case of the Stokes equations.

As for other Carleman estimates for the Navier-Stokes equations, we refer to
Choulli, Imanuvilov, Puel and Yamamoto \cite{ChImPuYa},
Fern\'andez-Cara, Guerrero, Imanuvilov and Puel \cite{FeGuImPu},
where the authors use a weight function in the form
$$
\exp\left( \frac{2sw(x)}{t(T-t)}\right)
$$
with some function $w$ and the weight function
decays to $0$ at $t=0,T$ exponentially.  Their Carleman
estimates hold over the whole domain $Q$ for $v$ satisfying
$v=0$ on $\ppp\Omega$ but not necessarily $v(\cdot,0) = v(\cdot,T) = 0$.
Those global Carleman estimate is convenient for proving
the Lipschitz stability for an inverse source problem
(e.g., \cite{ChImPuYa}) and the exact null controllability (\cite{FeGuImPu}),
but is not suitable for proving the unique continuation, and such a
weight function does not admit Carleman estimates for the Navier-Stokes
equations coupled with first-order equation or hyperbolic equation such as
a conservation law.
As for Carleman estimates for the Navier-Stokes equations, see also
Fan, Di Cristo, Jiang and Nakamura \cite{Fan1} and
Fan, Jiang and Nakamura \cite{Fan2} with extra data in a neighborhood
of the whole boundary, which is too much by considering the parabolicity of
the equations.

\section{Proof of Theorem 1}

{\bf First Step.}\\
Let $E \subset \R^n$ be a bounded domain with smooth boundary
$\ppp E$ and let
$E_{\delta} := \{ x\in E;\thinspace
\mbox{dist $(x, E) > \delta$}\}$ with small $\delta > 0$.

We prove\\
{\bf Lemma 1}.\\
{\it Let $p \in H^1(E)$ satisfy
$$
\Delta p = f_0 + \sum_{j=1}^n \ppp_jf_j \quad
\mbox{in $E$}
$$
and $\mbox{supp}\, p \subset E_{\delta}$.  Let $d_0 \in C^2(\overline{E})$
satisfy $d_0(x) > 0$ for $x \in E$ and $\vert \nabla d(x)\vert > 0$ for $x \in
E_{\delta}$.  We set $\va_0(x) = e^{\lambda d_0(x)}$ with large
constant $\lambda>0$.  Then there exist constants $C>0$ and
$s_1>0$ such that
$$
\int_E \left( \frac{1}{s}\vert \nabla p\vert^2
+ s\vert p\vert^2\right) e^{2s\va_0(x)} dx
\le C\int_E \left( \frac{1}{s^2}\vert f_0\vert^2
+ \sum_{j=1}^n \vert f_j\vert^2\right) e^{2s\va_0(x)} dx
$$
for all $s \ge s_1$.  The constants $C$ and $s_1$ are independent
of choices of $p$.}
\\
{\bf Proof.}  Since $\ppp E$ is of $C^3$-class, we choose
a function $\mu \in C^3(\overline{E})$ such that $0 \le \mu \le 1$,
$\mu > 0$ in $E$ and
$\mu = \left\{
\begin{array}{rl}
0, \quad &\mbox{in $\R^n \setminus E$},\\
1, \quad &\mbox{in $E_{\delta/2}$}.
\end{array}\right.$.
We set $\www{d_0}(x) = \mu(x) d_0(x)$ and $\www{\va_0}(x)
= e^{\lambda \www{d_0}(x)}$ for $x \in \overline{E}$.
Then $\www{d_0}(x) = 0$ for $x \in \ppp E$ and $\www{d_0}>0$,
$\vert \nabla\www{d_0}\vert = \vert \mu\nabla d_0 + d_0\nabla \mu\vert
= \vert \mu\nabla d_0\vert > 0$ in $E_{\delta}$.
Hence the $H^{-1}$-Carleman estimate for an elliptic operator
by Imanvilov and Puel \cite{ImPu} yields
$$
\int_E \left( \frac{1}{s}\vert \nabla p\vert^2
+ s\vert p\vert^2\right) e^{2s\www{\va_0}(x)} dx
\le C\int_E \left( \frac{1}{s^2}\vert f_0\vert^2
+ \sum_{j=1}^n \vert f_j\vert^2\right) e^{2s\www{\va_0}(x)} dx
$$
for all $s \ge s_1$.
Here we note that in Theorem 1.2 in \cite{ImPu}, we set $\omega = E\setminus
\overline{E_{\delta}}$ and use $p\vert_{\omega} = 0$.
Since $p=0$ in $E \setminus E_{\delta}$ and
$\www{d_0} = d_0$ in $E_{\delta}$, we complete the proof of Lemma 1.
\\
\vspace{0.2cm}
\\
{\bf Lemma 2.}\\
{\it There exist constants $s_0 > 0$ and $C>0$ such that
$$
\Vert (v,p)\Vert^2_{ \mathcal X_s(Q)}\le  C\int_Q \vert F\vert^2 \weight dxdt
+ C\int_Q (\vert h\vert^2 + \vert \nabla_{x,t}h\vert^2) \weight dxdt
                                                             \eqno{(2.1)}
$$
for all $s \ge s_0$ and $(v,p) \in H^{2,1}(Q) \times H^{1,0}(Q)$
satisfying (1.1),
\begin{align*}
&v(\cdot,0) = v(\cdot,T) = 0 \quad \mbox{in $\Omega$},\\
& \vert v\vert = \vert \nabla v\vert = \vert p\vert = 0
\quad \mbox{in $\ppp\Omega \times (0,T)$},
\end{align*}
and
$$
\ddd v = h \quad \mbox{in $Q$}
$$
with some $h \in H^{1,1}(Q)$.}
\\
{\bf Proof of Lemma 2.}
Thanks to the large parameter $s>0$, in view of (1.3),
it is sufficient to prove Lemma 1 for $B=0$ in (1.1).
In fact, the Carleman estimate with $B\ne 0$ follows from the case
of $B=0$ by replacing $F$ by $F-(v\cdot\nabla)B$
and estimating $\vert (F-(v\cdot\nabla)B)(x,t)\vert
\le \vert F(x,t)\vert + C\vert v(x,t)\vert$ for
$(x,t) \in Q$.  Then, choosing $s_0>0$ large, we can
absorb the term $\int_Q \vert v\vert^2 \weight dxdt$ into the
left-hand side of the Carleman estimate.

By the density argument, it is sufficient to prove the lemma for
$(v,p)$ such that $\mbox{supp} \,v$ and $\mbox{supp}\,p$ are compact in $Q$.
We consider
$$
\ppp_tv = \kappa\Delta v - (A\cdot\nabla)v - \nabla p + F
                                      \eqno{(2.2)}
$$
and
$$
\ddd v = h \qquad \mbox{in $Q$}.             \eqno{(2.3)}
$$
Taking the divergence of (2.2) and using (2.3), we obtain
$$
\Delta p = -\sum_{j,k=1}^n \{ \ppp_j((\ppp_kA_j)v_k)
- (\ppp_j\ppp_kA_j)v_k\}
+ \ddd F - \ppp_th - (A\cdot\nabla)h + \kappa\ddd(\nabla h)
\quad \mbox{in $Q$}.                    \eqno{(2.4)}
$$
Here we used
$$
\ddd ((A\cdot\nabla)v)
= \sum_{j,k=1}^n \ppp_k(A_j\ppp_jv_k)
= \sum_{j=1}^n A_j\ppp_j\left(\sum_{k=1}^n \ppp_kv_k\right)
+ \sum_{j,k=1}^n (\ppp_kA_j)\ppp_jv_k
$$
$$
= A\cdot \nabla(\ddd v)
+ \sum_{j,k=1}^n \{ \ppp_j((\ppp_kA_j)v_k)
- (\ppp_j\ppp_kA_j)v_k\}.                         \eqno{(2.5)}
$$
Moreover on the right-hand side of (2.4), the term
$\kappa\ddd (\nabla h)$ is not in $L^2(Q)$ because we assume only
$h \in H^{1,1}(Q)$.  Thus we cannot apply a usual Carleman estimate
requiring $\Delta p \in L^2(Q)$, and we need the $H^{-1}$-Carleman estimate.

By a usual density argument, we can assume that supp $p \subset Q$. 
By supp $p \subset Q$, fixing $t \in [0,T]$, we apply Lemma 1
to (2.4) and obtain
$$
\int_{\Omega} \left( \frac{1}{s} \vert \nabla p(x,t)\vert^2
+ s\vert p(x,t)\vert^2\right) e^{2s\va(x,t_0)} dx
$$
$$
\le C\int_{\Omega} (\vert F\vert^2 + \vert \ppp_th\vert^2
+ \vert \nabla h\vert^2 + \vert h\vert^2)e^{2s\va(x,t_0)} dx
+ C\int_{\Omega} \vert v(x,t)\vert^2 e^{2s\va(x,t_0)} dx
                                              \eqno{(2.6)}
$$
for $s \ge s_1$ where $s_1>0$ is a sufficiently large constant.

Let $s_0 := s_1e^{\lambda \beta T^2}$.  Then, $s \ge s_0$ implies
$$
se^{-\lambda\beta (t-t_0)^2} \ge se^{-\lambda\beta T^2} \ge s_1
$$
for $0 \le t \le T$, so that for fixed $t \in [0,T]$ by replacing
$s$ by $se^{-\la\beta (t-t_0)^2}$, by (2.5) we can see
\begin{align*}
&\int_{\Omega} \left( \frac{1}{s} \vert \nabla p(x,t)\vert^2
+ s\vert p(x,t)\vert^2\right) \exp( 2(se^{-\la\beta (t-t_0)^2})
\va(x,t_0)) dx \\
\le &C\int_{\Omega} (\vert F\vert^2 + \vert \ppp_th\vert^2
+ \vert \nabla h\vert^2 + \vert h\vert^2)\exp( 2(se^{-\la\beta (t-t_0)^2})
\va(x,t_0)) dx\\
+ & C\int_{\Omega} \vert v(x,t)\vert^2
\exp( 2(se^{-\la\beta (t-t_0)^2})\va(x,t_0))dx,
\end{align*}
that is,
\begin{align*}
& \int_{\Omega} \left( \frac{1}{s} \vert \nabla p(x,t)\vert^2
+ s\vert p(x,t)\vert^2\right) e^{2s\va(x,t)} dx \\
\le& C\int_{\Omega} (\vert F\vert^2 + \vert \ppp_th\vert^2
+ \vert \nabla h\vert^2 + \vert h\vert^2)e^{2s\va(x,t)} dx
+ C\int_{\Omega} \vert v(x,t)\vert^2 e^{2s\va(x,t)} dx
\end{align*}
for $s \ge s_0$ and $0 \le t \le T$.  Integrating this
inequality in $t$ over $(0,T)$,
we have
$$
\int_Q \left( \frac{1}{s} \vert \nabla p\vert^2
+ s\vert p\vert^2\right) e^{2s\va} dxdt
$$
$$
\le C\int_Q (\vert F\vert^2 + \vert \ppp_th\vert^2
+ \vert \nabla h\vert^2 + \vert h\vert^2)e^{2s\va} dxdt
+ C\int_Q \vert v\vert^2 e^{2s\va} dxdt           \eqno{(2.7)}
$$
for all $s \ge s_0$.

Next, regarding $F-\nabla p$ in (2.2) as non-homogeneous term,
we apply a Carleman estimate for the parabolic operator
$\ppp_tv - \kappa\Delta v + (A\cdot\nabla)v$ (e.g.,
Theorem 3.1 in Yamamoto \cite{Ya}) to (2.2):
$$
\frac{1}{s}\int_Q \left\{ \frac{1}{s}\left( \vert \ppp_tv\vert^2
+ \sum_{i,j=1}^n \vert \ppp_i\ppp_jv\vert^2\right)
+ s\vert \nabla v\vert^2 + s^3\vert v\vert^2\right\} \weight dxdt
$$
$$
\le C\int_Q \frac{1}{s}\vert \nabla p\vert^2 \weight dxdt
+ C\int_Q \frac{1}{s} \vert F\vert^2 \weight dxdt.
                                               \eqno{(2.8)}
$$

Substituting (2.7) into (2.8), we obtain
\begin{align*}
& \int_Q \left\{ \frac{1}{s^2}\left( \vert \ppp_tv\vert^2
+ \sum_{i,j=1}^n \vert \ppp_i\ppp_jv\vert^2\right)
+ \vert \nabla v\vert^2 + s^2\vert v\vert^2\right\} \weight dxdt \\
\le & C\int_Q \vert F\vert^2 \weight dxdt
+ C\int_Q (\vert \ppp_th\vert^2 + \vert \nabla h\vert^2
+ \vert h\vert^2)\weight dxdt\\
+ &C\int_Q \vert v\vert^2 \weight dxdt
+ \frac{C}{s}\int_Q \vert F\vert^2 \weight dxdt.
\end{align*}
Choosing $s_0 > 0$ large, we can absorb the third term on the right-hand
side into the left-hand side, again with (2.7), we complete the proof of
Lemma 2. $\blacksquare$
\\
\vspace{0.2cm}
\\
{\bf Second Step.}\\
Without loss of generality, we can assume that $d>0$ in $\Omega$ because
we replace $d$ by $d+C_0$ with large constant $C_0 > 0$ if necessary.

In this step, we will prove
\\
{\bf Lemma 3}.\\
{\it There exist constants $s_0>0$ and $C>0$ such that
\begin{align*}
& \Vert (v,p)\Vert^2_{ \mathcal X_s(D)}\\
\le& C\int_D \vert F\vert^2 \weight dxdt
+ C\int_D (\vert h\vert^2 + \vert \nabla_{x,t}h\vert^2)\weight dxdt\\
+ &Ce^{Cs}(\Vert v\Vert^2_{L^2(0,T;H^{\frac{3}{2}}(\Gamma))}
+ \Vert \ppp_tv\Vert^2_{L^2(0,T;H^{\frac{1}{2}}(\Gamma))}
+ \Vert \ppp_{\nu}v \Vert^2_{L^2(0,T;H^{\frac{1}{2}}(\Gamma))}
+ \Vert p\Vert^2_{L^2(0,T;H^{\frac{1}{2}}(\Gamma))})
\end{align*}
for all $s\ge s_0$ and $(v,p)\in H^{2,1}(D)\times H^{1,0}(D)$
satisfying (1.1), (1.6) and
$$
\ddd v = h \quad \mbox{in $D$}.                \eqno{(2.9)}
$$}
\\
{\bf Proof of Lemma 3.}
We take the zero extensions of $v,p, A, F$ to $Q$ from $D$ and
by the same letters we denote them:
$$
v
= \left\{
\begin{array}{rl}
v \quad & \mbox{on $\overline{D}$},\\
0 \quad & \mbox{in $Q \setminus D$},
\end{array}
\right.
\quad
p
= \left\{
\begin{array}{rl}
p \quad & \mbox{on $\overline{D}$},\\
0 \quad & \mbox{in $Q \setminus D$}, \thinspace
\mbox{etc.}\\
\end{array}
\right.
$$
By (1.6) we easily see that
$$
\ppp_i v
= \left\{
\begin{array}{rl}
\ppp_i v \quad & \mbox{on $\overline{D}$},\\
0 \quad & \mbox{in $Q\setminus D$},
\end{array}
\right.
\quad
\ppp_t v
= \left\{
\begin{array}{rl}
\ppp_t v \quad & \mbox{on $\overline{D}$},\\
0 \quad & \mbox{in $Q\setminus D$},
\end{array}
\right.
\quad
\ppp_i\ppp_j v
= \left\{
\begin{array}{rl}
\ppp_i\ppp_j v, \quad & \mbox{on $\overline{D}$},\\
0, \quad & \mbox{in $Q \setminus D$},
\end{array}
\right.
$$
and
$$
\ppp_i p
= \left\{
\begin{array}{rl}
\ppp_i p, \quad & \mbox{on $\overline{D}$},\\
0, \quad & \mbox{in $Q\setminus D$}
\end{array}
\right.
$$
for $1 \le i,j \le n$.
Moreover, since $v=0$ on $\ppp D \setminus
(\Gamma \times (0,T))$ by (1.6), setting
\\
$h = \left\{
\begin{array}{rl}
h \quad & \mbox{on $\overline{D}$},\\
0 \quad & \mbox{in $Q\setminus D$}
\end{array}
\right.
$, we see that $h \in H^{1,1}(Q)$ and
$$
\ddd v = h \quad \mbox{in $Q$}    \eqno{(2.10)}
$$
and
$$
\ppp_tv = \kappa\Delta v + (A\cdot \nabla)v + \nabla p + F
\quad \mbox{in $Q$}.                                \eqno{(2.11)}
$$

By the Sobolev extension theorem, there exist $\www{p}
\in L^2(0,T;H^1(\Omega))$ and $v \in H^{2,1}(Q)$ such that
$$
\left\{
\begin{array}{rl}
&\www{v} = v, \thinspace \ppp_{\nu}\www{v} = \ppp_{\nu} v, \thinspace
\www{p} = p \quad \mbox{on $\ppp\Omega \times (0,T)$},\\
& \mbox{supp $\www{v}(x,\cdot) \subset (0,T)$ for almost all $x \in \Omega$}
\end{array}
\right.                      \eqno{(2.12)}
$$
and
$$
\Vert \www{v}\Vert_{H^{2,1}(Q)}
+ \Vert \ppp_t{\www{v}}\Vert_{L^2(0,T;H^1(\Omega))}
+ \Vert \www{p}\Vert_{L^2(0,T;H^1(\Omega))}
$$
$$
\le C(\Vert v\Vert_{L^2(0,T;H^{\frac{3}{2}}(\Gamma))}
+ \Vert \ppp_t v\Vert_{L^2(0,T;H^{\frac{1}{2}}(\Gamma))}
+ \Vert \ppp_{\nu} v\Vert_{L^2(0,T;H^{\frac{1}{2}}(\Gamma))}
+ \Vert p\Vert_{L^2(0,T;H^{\frac{1}{2}}(\Gamma))}).      \eqno{(2.13)}
$$
The last condition in (2.12) can be seen by $v(\cdot, 0) =
v(\cdot,T) = 0$ in $\Omega$ which follows from (1.6).

We set
$$
u = v - \www{v}, \quad q = p-\www{p} \quad \mbox{in $Q$}.
$$
Then, in view of (2.10) - (2.12), we have
$$
\vert u\vert  = \vert \nabla u\vert = \vert q\vert = 0 \quad
\mbox{on $\ppp\Omega\times (0,T)$}          \eqno{(2.14)}
$$
and
$$
\ppp_tu - \kappa\Delta u + \nabla q + (A\cdot\nabla)u
= F - (\ppp_t\www{v} - \kappa\Delta \www{v} + (A\cdot\nabla)\www{v}
+ \nabla\www{p}) =: G \quad \mbox{in $Q$},    \eqno{(2.15)}
$$
$$
\ddd u = h -\ddd \www{v} \in H^{1,1}(Q).              \eqno{(2.16)}
$$

We choose a bounded domain $\www{\Omega}$ with smooth boundary
$\ppp\www{\Omega}$ such that $\www{\Omega} \supset \Omega$,
$\overline{\Gamma} = \ppp\Omega \cap \www{\Omega}$ and
$\ppp\widetilde{\Omega} \cap \overline{\Omega}
= \ppp\Omega \setminus \Gamma$.  In other words, the
domain $\widetilde{\Omega}$ is constructed by expanding $\Omega$ only over
$\Gamma$ to the exterior such that the boundary $\ppp\widetilde{\Omega}$
is smooth.  We set
$$
\www{Q} =\www{\Omega} \times (0,T).
$$

Let us recall that $d$ satisfies (1.4).  Since we can further
choose $\widetilde{\Omega}$ such that
$\widetilde{\Omega}\setminus \Omega$ is included in a sufficiently small
ball, we see that there exists
an extension $\widetilde{d}$ in $\widetilde{\Omega}$ of $d$
satisfying $\vert \nabla\www{d}\vert > 0$ in $\www{\Omega}$.

We take the zero extensions of $u, q, A, G$ and $h-\ddd \www{v}$
to $\www{\Omega}$ and by the same letters we denote them.
Therefore by (2.14) - (2.16), the zero extensions of $u$ and $h - \ddd\www{v}$
satisfies
$$
\ddd u = h - \ddd\widetilde{v} \in H^{1,1}(\www{Q})    \eqno{(2.17)}
$$
and
$$
\ppp_tu - \kappa\Delta u + \nabla q + (A\cdot\nabla)u = G \quad
\mbox{in $\www{Q}$}.                               \eqno{(2.18)}
$$
By the zero extensions and (1.6), we obtain
$$
u(\cdot,0) = u(\cdot, T) = 0 \quad \mbox{in $\www{\Omega}$},
$$
$$
\vert u\vert = \vert \nabla u\vert = \vert q\vert = 0 \quad
\mbox{on $\ppp\www{\Omega} \times (0,T)$}.        \eqno{(2.19)}
$$
Therefore, by noting (2.19), we apply Lemma 2 to (2.17)
and (2.18), and we obtain
\begin{align*}
& \Vert (u,q)\Vert^2_{ \mathcal X_s(\tilde Q)}
\le  C\int_{\www{Q}} \vert G\vert^2 \weight dxdt\\
+ &C\int_{\www{Q}} (\vert h - \ddd\www{v}\vert^2
+ \vert \nabla_{x,t} (h-\ddd\www{v})\vert^2)\weight dxdt
\end{align*}
for $s \ge s_0$.
Hence
\begin{align*}
& \Vert (v-\tilde v,p-\tilde p)\Vert^2_{ \mathcal X_s(Q)} \\
\le & C\int_Q \vert F\vert^2 \weight dxdt\\
+ &C\int_Q \vert \ppp_t\www{v} - \kappa\Delta\www{v}
+ (A\cdot\nabla)\www{v} + \nabla\www{p}\vert^2 \weight dxdt\\
+& C\int_Q \left( \vert h\vert^2 + \vert \nabla_{x,t}h\vert^2
+ \vert \nabla\www{v}\vert^2
+ \sum_{i,j=1}^n \vert \ppp_i\ppp_j\www{v}\vert^2
+ \vert \nabla(\ppp_t\www{v})\vert^2\right)\weight dxdt.
\end{align*}
Using $\vert \ppp_tv\vert^2 \le 2\vert \ppp_t\www{v}\vert^2
+ 2\vert \ppp_t(v-\www{v})\vert^2$, etc. on the left-hand side, we
have
\begin{align*}
& \Vert (v,p)\Vert^2_{ \mathcal X_s(Q)}
\le 2\Vert (\tilde v,\tilde p)\Vert^2_{ \mathcal X_s(Q)} \\
+ &2C\int_Q \vert F\vert^2 \weight dxdt\\
+ &2C\int_Q \vert \ppp_t\www{v} - \kappa\Delta\www{v}
+ (A\cdot\nabla)\www{v} + \nabla\www{p}\vert^2 \weight dxdt\\
+& 2C\int_Q \left( \vert h\vert^2 + \vert \nabla_{x,t}h\vert^2
+ \vert \nabla\www{v}\vert^2
+ \sum_{i,j=1}^n \vert \ppp_i\ppp_j\www{v}\vert^2
+ \vert \nabla(\ppp_t\www{v})\vert^2\right)\weight dxdt\\
\le& C\int_Q \vert F\vert^2 \weight dxdt\\
+& Ce^{Cs}(\Vert \www{v}\Vert^2_{H^{2,1}(Q)}
+ \Vert \www{p}\Vert^2_{H^{1,0}(Q)}
+ \Vert \nabla\ppp_t\www{v}\Vert^2_{L^2(Q)})\\
+& C\int_Q (\vert h\vert^2 + \vert \nabla_{x,t}h\vert^2) \weight dxdt
\end{align*}
for $s \ge s_0$.  Since $F$ and $h$ are zero outside of $D$,
in view of (2.13), the proof of Lemma 3 is completed. $\blacksquare$
\\
\vspace{0.2cm}
\\

{\bf Third Step.}\\
For $r>0$ and $x_0 \in \R^n$, we set $B_r(x_0):= \{x \in \R^n; \thinspace
\vert x-x_0\vert < r\}$.
Then we prove
\\
{\bf Lemma 4.}\\
{\it Let $v \in H^2(\Omega)$ and $p \in H^1(\Omega)$.\\
{\bf (1) Case $n=3$}:  For any $x_0 \in \ppp\Omega$, there exist
$r>0$ and a $10\times 10$ matrix $A \in C^1(\overline{B_r(x_0)})$ such that
$$
\ppp\Omega \cap B_r(x_0) = \{ x(\theta_1,\theta_2);
\thinspace (\theta_1,\theta_2) \in D_1\}
$$
where $x(\theta_1,\theta_2) = (x_1(\theta_1,\theta_2),
x_2(\theta_1,\theta_2), x_3(\theta_1,\theta_2)) \in \R^3$,
$D_1 \subset \R^2$ is a bounded domain and the functions $x_1, x_2, x_3$ 
with respect to $\theta_1, \theta_2$ are in $C^2(\overline{D_1})$ and
$$
\mbox{det}\thinspace A(x(\theta_1,\theta_2)) \ne 0, \quad
(\theta_1, \theta_2) \in \overline{D_1}
$$
and
$$
A(x(\theta_1,\theta_2))
\left(
\begin{array}{cccc}
&(\nabla_xv_1)(x(\theta_1,\theta_2))\\\
&(\nabla_xv_2)(x(\theta_1,\theta_2))\\
&(\nabla_xv_3)(x(\theta_1,\theta_2))\\
&p(x(\theta_1,\theta_2))\\
\end{array}\right)
= \left(
\begin{array}{ccccc}
&\nabla_{\theta_1,\theta_2}(v_1(x(\theta_1,\theta_2)))\\
&\nabla_{\theta_1,\theta_2}(v_2(x(\theta_1,\theta_2)))\\
&\nabla_{\theta_1,\theta_2}(v_3(x(\theta_1,\theta_2)))\\
&(\sigma(v,p)\nu)(x(\theta_1,\theta_2))\\
&(\ddd v)(x(\theta_1,\theta_2))\\
\end{array}\right),
\quad (\theta_1, \theta_2) \in D_1.
$$
\\
{\bf (2) Case $n=2$}:  For any $x_0 \in \ppp\Omega$, there exist
$r>0$ and a $5\times 5$ matrix $A \in C^1(\overline{B_r(x_0)})$ such that
$$
\ppp\Omega \cap B_r(x_0) = \{ x(\theta_1); \thinspace
\theta_1 \in I_1\}
$$
where $x(\theta_1) := (x_1(\theta_1), x_2(\theta_1)) \in \R^2$,
$I_1 \subset \R$ is an open interval, and the functions $x_1, x_2$ are in
$C^2(\overline{I_1})$, and
$$
\mbox{det}\thinspace A(x(\theta_1)) \ne 0, \quad \theta_1 \in \overline{I_1}
$$
and
$$
A(x(\theta_1))
\left(
\begin{array}{ccc}
&(\nabla_xv_1)(x(\theta_1))\\
&(\nabla_xv_2)(x(\theta_1))\\
&p(x(\theta_1))\\
\end{array}\right)
= \left(
\begin{array}{cccc}
&\frac{d}{d\theta_1}v_1(x(\theta_1))\\
&\frac{d}{d\theta_1}v_2(x(\theta_1))\\
&(\sigma(v,p)\nu)(x(\theta_1))\\
&(\ddd v)(x(\theta_1))\\
\end{array}\right),
\quad \theta_1 \in I_1.
$$}
\\
\vspace{0.2cm}
{\bf Remark.}
The lemma guarantees that the boundary data $(v,\ppp_{\nu}v,p)$ and
$(v,\sigma(v,p)\nu)$ are equivalent (e.g., Imanuvilov and Yamamoto
\cite{ImYa2}).  As related papers on inverse boundary value problems
for the Navier-Stokes equations in view of this equivalence, see
Imanuvilov and Yamamoto \cite{ImYa1}, Lai, Uhlmann and
Wang \cite{LUW}.
\\
{\bf Proof of Lemma 4.}
We prove only in the case of $n=3$.  The case of $n=2$ is similar and simpler.
It is sufficient to consider only on
a sufficiently small subboundary $\Gamma_0$ of
$\ppp\Omega$.  Without loss of generality, we can
assume that $\Gamma_0$ is represented by by $(x_1,x_2,\gamma(x_1,x_2))$ where
$\gamma \in C^2(\overline{D_1})$,
$\theta_1 = x_1$, $\theta_2 = x_2$,
$x_3 = \gamma(x_1,x_2)$ for $(x_1,x_2) \in D_1$.  Moreover we assume that
$\Omega$ is located upper $x_3 = \gamma(x_1,x_2)$.

By the density argument, we can assume that $v \in C^1(\overline{\Omega})$ and
$p \in C(\overline{\Omega})$.

We set $\gamma_1;= \ppp_1\gamma$ and $\gamma_2:=
\ppp_2\gamma$.
On $\Gamma_0$, we have
$$
\nu(x) = \frac{1}{1+\gamma_1^2+\gamma_2^2}
\left(
\begin{array}{ccc}
&\gamma_1\\
&\gamma_2\\
&-1\\
\end{array}\right).                 \eqno{(2.20)}
$$

By the definition, we have
$$
\sigma(v,p)\nu = \kappa
\left(
\begin{array}{ccc}
2\ppp_1v_1 - \frac{p}{\kappa} & \ppp_1v_2+\ppp_2v_1 &
\ppp_1v_3+\ppp_3v_1 \\
\ppp_1v_2+\ppp_2v_1 & 2\ppp_2v_2-\frac{p}{\kappa} &
\ppp_2v_3+\ppp_3v_2 \\
\ppp_1v_3 + \ppp_3v_1 & \ppp_2v_3+\ppp_3v_2 &
2\ppp_3v_3-\frac{p}{\kappa} \\
\end{array}\right)
\left(
\begin{array}{ccc}
& \nu_1\\
& \nu_2\\
& \nu_3\\
\end{array}\right)                           \eqno{(2.21)}
$$
$$
=:
\left(
\begin{array}{ccc}
& q_1\\
& q_2\\
& q_3\\
\end{array}\right).
$$
We further set
$$
q_4 := (\ddd v)(x_1,x_2,\gamma(x_1,x_2)),
$$
$$
g_k(x_1,x_2) := v_k\XXX, \quad k=1,2,3.
$$
Then
$$
\ppp_1g_k = \ppp_1v_k + \gamma_1(\ppp_3v_k)\XXX,
$$
$$
\ppp_2g_k = \ppp_2v_k + \gamma_2(\ppp_3v_k)\XXX,
$$
that is,
$$
\left\{
\begin{array}{rl}
&\ppp_1v_k\XXX = \ppp_1g_k - \gamma_1(\ppp_3v_k)\XXX, \\
&\ppp_2v_k\XXX = \ppp_2g_k - \gamma_2(\ppp_3v_k)\XXX, \quad k=1,2,3,\\
\end{array}\right.
                                     \eqno{(2.22)}
$$
and
$$
(\ppp_3v_3)\XXX = q_4 - (\ppp_1v_1+\ppp_2v_2)\XXX     \eqno{(2.23)}
$$
for $(x_1,x_2) \in D_1$.
Setting
$$
\left\{
\begin{array}{rl}
& h_1(x_1,x_2) = (\ppp_3v_1)\XXX, \\
& h_2(x_1,x_2) = (\ppp_3v_2)\XXX, \\
\end{array}\right.
                                     \eqno{(2.24)}
$$
by (2.22) and (2.23) we obtain
\begin{align*}
& (\ppp_3v_3)\XXX
= q_4 - (\ppp_1v_1+\ppp_2v_2)\XXX\\
=& (q_4 - \ppp_1g_1 - \ppp_2g_2)(x_1,x_2)
+ (\gamma_1\ppp_3v_1+\gamma_2\ppp_3v_2)\XXX\\
=: & g_0(x_1,x_2) + (\gamma_1\ppp_3v_1+\gamma_2\ppp_3v_2)\XXX
\end{align*}
$$
= (g_0+\gamma_1h_1+\gamma_2h_2)(x_1,x_2)    \eqno{(2.25)}
$$
and so
$$
\left\{
\begin{array}{rl}
& (\ppp_1v_1)\XXX = (\ppp_1g_1 - \gamma_1h_1)\XXX,\\
& (\ppp_2v_1)\XXX = (\ppp_2g_1 - \gamma_2h_1)\XXX,\\
& (\ppp_1v_2)\XXX = (\ppp_1g_2 - \gamma_1h_2)\XXX,\\
& (\ppp_2v_2)\XXX = (\ppp_2g_2 - \gamma_2h_2)\XXX,\\
& (\ppp_1v_3)\XXX = \ppp_1g_3 - \gamma_1g_0 - \gamma_1^2h_1
- \gamma_1\gamma_2h_2,\\
& (\ppp_2v_3)\XXX = \ppp_2g_3 - \gamma_2g_0 - \gamma_1\gamma_2h_1
- \gamma_2^2h_2,\quad (x_1,x_2) \in D_1.
\end{array}\right.
\eqno{(2.26)}
$$
On the other hand, (2.21) yields
$$
\frac{1+\gamma_1^2+\gamma_2^2}{\kappa}q_1
= (2\gamma_1\ppp_1v_1+\gamma_2\ppp_1v_2+\gamma_2\ppp_2v_1
- \ppp_1v_3 - \ppp_3v_1)\XXX - \frac{\gamma_1}{\kappa}p,
$$
$$
\frac{1+\gamma_1^2+\gamma_2^2}{\kappa}q_2
= (\gamma_1\ppp_1v_2+\gamma_2\ppp_2v_1 + 2\gamma_2\ppp_2v_2
- \ppp_2v_3 - \ppp_3v_2)\XXX - \frac{\gamma_2}{\kappa}p
$$
and
$$
\frac{1+\gamma_1^2+\gamma_2^2}{\kappa}q_3
= (\gamma_1\ppp_1v_3+\gamma_1\ppp_3v_1 + \gamma_2\ppp_2v_3
+ \gamma_2\ppp_3v_2 - 2\ppp_3v_3)\XXX + \frac{1}{\kappa}p,
\quad (x_1,x_2) \in D_1.
$$
Substitute (2.25) and (2.26), we have
$$
\left\{
\begin{array}{rl}
& \frac{1+\gamma_1^2+\gamma_2^2}{\kappa}q_1
= -(1+\gamma_1^2+\gamma_2^2)h_1 - \frac{\gamma_1}{\kappa}p+ G_1, \\
& \frac{1+\gamma_1^2+\gamma_2^2}{\kappa}q_2
= -(1+\gamma_1^2+\gamma_2^2)h_2 - \frac{\gamma_2}{\kappa}p+ G_2, \\
& \frac{1+\gamma_1^2+\gamma_2^2}{\kappa}q_3
= -\gamma_1(1+\gamma_1^2+\gamma_2^2)h_1
  -\gamma_2(1+\gamma_1^2+\gamma_2^2)h_2 + \frac{1}{\kappa}p + G_3.
\end{array}
\right.                                  \eqno{(2.27)}
$$
Here $G_k$, $k=1,2,3$, are linear combinations of $\ppp_jg_k, q_1,q_2,
q_3, q_4$, $j=1,2, k=1,2,3$, with coefficients given by $\gamma$ and
its first-order derivatives.
We can uniquely solve (2.27) with respect to $h_1,h_2, p$:
$$
\left(
\begin{array}{ccc}
& h_1(x_1,x_2)\\
& h_2(x_2,x_2)\\
& p\XXX\\
\end{array}\right)
= \widetilde{A}(x_1,x_2)
\left(
\begin{array}{ccc}
& \frac{1+\gamma_1^2+\gamma_2^2}{\kappa}q_1 - G_1\\
& \frac{1+\gamma_1^2+\gamma_2^2}{\kappa}q_2 - G_2\\
& \frac{1+\gamma_1^2+\gamma_2^2}{\kappa}q_3 - G_3\\
\end{array}\right), \quad (x_1,x_2) \in D_1.
                                    \eqno{(2.28)}
$$
Here $\widetilde{A} \in C^1(\overline{D_1})$ and
 $\mbox{det}\,\widetilde{A} \ne 0$ on $\overline{D_1}$.  The equations
(2.25), (2.26) and (2.28) imply the existence of a $10\times 10$ matrix
$A \in C^1(\overline{D_1})$ satisfying the conditions in the lemma.
Thus the proof of Lemma 4 is completed. $\blacksquare$
\\
\vspace{0.2cm}

Now, in terms of Lemmata 3 and 4, we complete the proof of Theorem 1
as follows.
We consider only the case of $n=3$.  Without loss of generality,
$\Gamma$ is given by
$\Gamma = \{ \XXX;\thinspace x_1,x_2 \in D_1\}$ with
$\gamma \in C^2(\overline{D_1})$.
\\
We set $\nabla_{x_1,x_2}v = (\ppp_1v_1,\ppp_2v_1,
\ppp_1v_2, \ppp_2v_2, \ppp_1v_3, \ppp_2v_3)^T$.
Then, by Lemmata 2 and 3, we have
\begin{align*}
& \left(
\begin{array}{cc}
& \ppp_{\nu}v\XXX\\
& p\XXX\\
\end{array}\right)
= \left(
\begin{array}{cc}
& \frac{1}{1+\gamma_1^2+\gamma_2^2}
((\ppp_1\gamma)\ppp_1v + (\ppp_2\gamma)\ppp_2v - \ppp_3v)\XXX\\
& p\XXX\\
\end{array}\right)\\
=& B_1(x_1,x_2)
\left(
\begin{array}{ccc}
&(\nabla_{x_1,x_2}v)\XXX\\
&(\sigma(v,p)\nu)\XXX\\
\end{array}\right), \quad (x_1,x_2) \in D_1,
\end{align*}
with a $4\times 6$ matrix $B_1 \in C^1(\overline{D_1})$.
Therefore
$$
\Vert \ppp_{\nu}v(\cdot,t)\Vert_{H^1(\Gamma)}
+ \Vert p(\cdot,t)\Vert_{H^1(\Gamma)}
= \left\Vert
B\left(
\begin{array}{ccc}
&\nabla_{x_1,x_2}v\\
&\sigma(v,p)\nu\\
\end{array}\right)
(\cdot,t)\right\Vert_{H^1(\Gamma)}
\le C\left\Vert
\left(
\begin{array}{ccc}
&\nabla_{x_1,x_2}v\\
&\sigma(v,p)\nu\\
\end{array}\right)
(\cdot,t)\right\Vert_{H^1(\Gamma)}
$$
and
$$
\Vert \ppp_{\nu}v(\cdot,t)\Vert_{L^2(\Gamma)}
\le C\left\Vert
\left(
\begin{array}{ccc}
&\nabla_{x_1,x_2}v\\
&\sigma(v,p)\nu\\
\end{array}\right)
(\cdot,t)\right\Vert_{L^2(\Gamma)}
$$
by $B \in C^1(\overline{D_1})$.  Consequently the interpolation
inequality (e.g., Theorem 7.7 (p.36) in Lions and Magenes \cite{LM})
yields
$$
\Vert \ppp_{\nu}v(\cdot,t)\Vert_{H^{\frac{1}{2}}(\Gamma)}
+ \Vert p(\cdot,t)\Vert_{H^{\frac{1}{2}}(\Gamma)}
\le \left\Vert
\left(
\begin{array}{ccc}
&\nabla_{x_1,x_2}v\\
&\sigma(v,p)\nu\\
\end{array}\right)
(\cdot,t)\right\Vert_{H^{\frac{1}{2}}(\Gamma)}
$$
for $0\le t \le T$.
Hence
$$
\Vert \ppp_{\nu}v\Vert_{L^2(0,T;H^{\frac{1}{2}}(\Gamma))}
+ \Vert p(\cdot,t)\Vert_{L^2(0,T;H^{\frac{1}{2}}(\Gamma))}
\le C(\Vert v\Vert_{L^2(0,T;H^1(\Gamma))}
+ \Vert \sigma(v,p)\nu\Vert_{L^2(0,T;H^{\frac{1}{2}}(\Gamma))}).
$$
With this, Lemma 3 completes the proof of Theorem 1. $\blacksquare$
\\
\section{Conditional stability for the lateral Cauchy problem}

In this section, we discuss \\
{\bf lateral Cauchy problem}\\
We are given a suboundary $\Gamma$ of $\ppp\Omega$ arbitrarily.
Let $(v,p) \in H^{2,1}(Q) \times H^{1,0}(Q)$ satisfy (1.1) and
(1.2).  Determine $(v,p)$ in some subdomain of $Q$ by
$(v, \sigma(v,p)\nu)$ on $\Gamma \times (0,T)$.

In the case of the parabolic equation, there are very many works, and
here we do not list up comprehensively and as restricted references,
see Landis \cite{La}, Mizohata \cite{Mi},
Saut and Scheurer \cite{SS}, Sogge \cite{So}.  See also the monographs
Beilina and Klibanov \cite{BeKl}, Isakov \cite{Is},
Klibanov and Timonov \cite{KlTi}.

Combining a Carleman estimate and a cut-off function, we can 
prove
\\
{\bf Proposition 1.}\\
{\it Let $\va(x,t)$ be given in Theorem 1.  We set
$$
Q(\ep) = \{ (x,t) \in \Omega \times (0,T); \va(x,t) >
\ep\}
$$
with $\ep > 0$.  Moreover we assume that
$$
\oooo{Q(0)} \subset Q \cup (\Gamma \times [0,T])
$$
with subboundary $\Gamma \subset \ppp\Omega$.  Then for any small
$\ep>0$, there exist constants $C>0$ and $\theta\in (0,1)$ such that
$$
\Vert v\Vert_{H^{2,1}(Q(\ep))} + \Vert p\Vert_{H^{1,0}(Q(\ep))}
\le C(\Vert v\Vert^{1-\theta}_{H^{1,1}(Q)}
+ \Vert p\Vert_{L^2(Q)})G^\theta
+ CG,
$$
where we set
$$
G^2 := \Vert F\Vert_{L^2(Q)}^2
+ \Vert v\Vert^2_{L^2(0,T;H^{\frac{3}{2}}(\Gamma))}
+ \Vert \ppp_tv \Vert^2_{L^2(0,T;H^{\frac{1}{2}}(\Gamma))}
+ \Vert \sigma(v,p)\nu\Vert^2_{L^2(0,T;H^{\frac{1}{2}}(\Gamma))}.
$$}

As for the proof of Proposition 1, see Theorem 3.2.2 in section 3.2 of
\cite{Is} for example.

Proposition 1 gives an estimate of the solution in $Q(\ep)$ by
data on $\Gamma \times (0,T)$, and $Q(\ep)$ and $\Gamma$ are
determined by an a priori given function $d(x)$.
Therefore the proposition does not give a suitable answer
to our lateral Cauchy problem as stated above, where we are
requested to estimate the solution by data on as a small subboundary
$\Gamma \times (0,T)$ as possible.

In fact, in this section, we prove
\\
{\bf Theorem 2 (conditional stability)}\\
{\it Let $\Gamma \subset \ppp\Omega$ be an arbitrary non-empty
subboundary of $\ppp\Omega$.
For any $\ep > 0$ and
an arbitrary bounded domain $\Omega_0$ such that
$\overline{\Omega_0} \subset \Omega \cup \Gamma$, $\ppp\Omega_0
\cap \ppp\Omega$ is a non-empty open subset of $\ppp\Omega$ and
$\ppp\Omega_0 \cap \ppp\Omega \subsetneqq \Gamma$,
there exist constants $C>0$ and $\theta \in (0,1)$ such that
$$
\Vert v\Vert_{H^{2,1}(\Omega_0\times (\ep,T-\ep))}
+ \Vert p\Vert_{H^{1,0}(\Omega_0\times (\ep,T-\ep))}
$$
$$
\le C(\Vert v\Vert_{H^{1,1}(Q)} + \Vert p\Vert_{L^2(Q)})
^{1-\theta}
(\Vert F\Vert_{L^2(Q)}
+ \Vert v\Vert_{L^2(0,T;H^{\frac{3}{2}}(\Gamma))}
+ \Vert v \Vert_{H^1(0,T;H^{\frac{1}{2}}(\Gamma))}
+ \Vert \sigma(v,p)\nu\Vert_{L^2(0,T;H^{\frac{1}{2}}(\Gamma))})^{\theta}
$$
$$
+ C(\Vert F\Vert_{L^2(Q)}
+ \Vert v\Vert_{L^2(0,T;H^{\frac{3}{2}}(\Gamma))}
+ \Vert \ppp_tv \Vert_{L^2(0,T;H^{\frac{1}{2}}(\Gamma))}
+ \Vert \sigma(v,p)\nu\Vert_{L^2(0,T;H^{\frac{1}{2}}(\Gamma))}).
                                                            \eqno{(3.1)}
$$}
\\

In Theorem 2, in order to estimate $(v,p)$, we have to assume
a priori bounds of $\Vert v\Vert_{H^{1,1}(Q)}$ and
$\Vert p\Vert_{L^2(Q)}$.  Thus estimate (3.1)
is called a conditional stability estimate. %
We note that (3.1) is rewritten as
\begin{align*}
&\Vert v\Vert_{H^{2,1}(\Omega_0\times (\ep,T-\ep))}
+ \Vert p\Vert_{H^{1,0}(\Omega_0\times (\ep,T-\ep))}\\
=& O((\Vert F\Vert_{L^2(Q)}
+ \Vert v\Vert_{L^2(0,T;H^{\frac{3}{2}}(\Gamma))}
+ \Vert v \Vert_{H^1(0,T;H^{\frac{1}{2}}(\Gamma))}
+ \Vert \sigma(v,p)\nu\Vert_{L^2(0,T;H^{\frac{1}{2}}(\Gamma))})^{\theta})
\end{align*}
as $\Vert F\Vert_{L^2(Q)} + \Vert v\Vert_{L^2(0,T;H^{\frac{3}{2}}(\Gamma))}
+ \Vert v \Vert_{H^1(0,T;H^{\frac{1}{2}}(\Gamma))}
+ \Vert \sigma(v,p)\nu\Vert_{L^2(0,T;H^{\frac{1}{2}}(\Gamma))}
\longrightarrow 0$.
Thus the estimate indicates stability of H\"older type. 

For the homogeneous Stokes equations:
$$
\ppp_tv - \Delta v + \nabla p = 0, \quad
\ddd v = 0 \quad \mbox{in $Q$},
$$
Boulakia \cite{Bo} (Proposition 2) proved the conditional stability
in $\Omega_0 \times (\ep, T-\ep)$ on the basis of a Carleman estimate
in \cite{Bo}.  The norm of boundary data in \cite{Bo} is stronger than our
chosen norm.

The theorem does not directly give an estimate when $\Omega_0 = \Omega$, but
we can derive an estimate in $\Omega$ by an argument similar to
Theorem 5.2 in Yamamoto \cite{Ya} and we do not discuss details.
Boulakia \cite{Bo} (Theorem 1) established a conditional stability
estimate up to $\ppp\Omega$ by boundary or interior data.
The argument is based on the interior estimate in $\Omega_0
\times (\ep, T-\ep)$ and an
argument similar to Theorem 5.2 in \cite{Ya}.

Theorem 2 immediately implies the global uniqueness of
the solution:
\\
{\bf Corollary}.\\
{\it Let $\Gamma \subset \ppp\Omega$ be an arbitrarily fixed subboundary.
If $(v,p) \in H^{2,1}(Q)\times H^{1,0}(Q)$ satisfies (1.1) and (1.2), and
$v = \sigma(v,p)\nu = 0$ on $\Gamma \times (0,T)$, then
$\vert v\vert = \sigma(v,p)\nu = 0$ in $\Omega \times (0,T)$.}
\\
\vspace{0.2cm}
\\
{\bf Proof of Theorem 2}.
Once a relevant Carleman estimate for the Navier-Stokes equations is proved,
the proof is similar to Theorem 5.1 in \cite{Ya}.
Thus, according to $\Omega_0$ and
$\Gamma$, we have to choose a suitable weight function $\va$.
For this, we show\\
{\bf Lemma 5.}\\
{\it Let $\omega$ be an arbitrarily fixed subdomain of $\Omega$
such that $\overline{\omega} \subset \Omega$.  Then there exists
a function $d \in C^2(\oooo\Omega)$ such that
$$
d(x) > 0 \quad x \in \Omega, \quad d\vert_{\ppp\Omega}
= 0, \quad \vert \nabla d(x)\vert > 0, \quad
x\in \overline{\Omega\setminus\omega}.
$$}
\\
\vspace{0.1cm}
\\
For the proof, see Fursikov and Imanuvilov \cite{FuIm}, Imanuvilov \cite{Im},
Imanuvilov, Puel and Yamamoto \cite{ImPuYa}.

We choose a bounded domain $\Omega_1$
with smooth boundary such that
$$
\Omega \subsetneqq \Omega_1, \quad
\oooo{\Gamma} = \oooo{\ppp\Omega \cap \Omega_1}, \quad
\ppp\Omega \setminus \Gamma \subset \ppp\Omega_1,      \eqno{(3.2)}
$$
and $\Omega_1\setminus \oooo\Omega$ contains
some non-empty open set.
We note that $\Omega_1$ is constructed by taking a union of
$\Omega$ and a domain $\www{\Omega} \subset \R^n \setminus
\overline{\Omega}$ such that
$\www\Omega \cap \ppp\Omega = \Gamma$.
Choosing $\oooo\omega \subset \Omega_1 \setminus \oooo\Omega$,
and applying Lemma 5 to obtain $d \in C^2(\oooo\Omega_1)$ satisfying
$$
d(x) > 0, \quad x \in \Omega_1, \quad
d(x) = 0, \quad x \in \ppp\Omega_1, \quad
\vert \nabla d(x)\vert > 0, \quad x \in \oooo\Omega.
                                                            \eqno{(3.3)}
$$
Then, since $\overline{\Omega_0} \subset \Omega_1$, we can choose
sufficiently large $N > 1$ such that
$$
\{ x\in \Omega_1; \thinspace d(x) > \frac{4}{N}\Vert d\Vert
_{C(\oooo{\Omega_1})} \} \cap \oooo\Omega \supset \Omega_0.                                                                    \eqno{(3.4)}
$$
Moreover we choose sufficiently large $\beta > 0$ such that
$$
\beta\ep^2 < \Vert d\Vert_{C(\oooo{\Omega_1})} < 2\beta\ep^2.     \eqno{(3.5)}
$$
We arbitrarily fix $t_0 \in [\sqrt{2}\ep, T-\sqrt{2}\ep]$.  We set
$\va(x,t) = e^{\la\psi(x,t)}$ with fixed large parameter $\la>0$ and
$\psi(x,t) = d(x) - \beta (t-t_0)^2$,
$\mu_k = \exp\left( \la\left(\frac{k}{N}\Vert d\Vert_{C(\oooo{\Omega_1})}
- \frac{\beta\ep^2}{N}\right)\right)$, $k=1,2,3,4$, and
$D = \{ (x,t); \thinspace x \in\oooo\Omega, \quad \va(x,t) > \mu_1
\}$.

Then we can verify that
$$
\Omega_0 \times \left( t_0 - \frac{\ep}{\sqrt{N}},
t_0 + \frac{\ep}{\sqrt{N}} \right) \subset D
\subset \oooo\Omega \times (t_0-\sqrt{2}\ep,t_0+\sqrt{2}\ep).
                                      \eqno{(3.6)}
$$
In fact, let $(x,t) \in\Omega_0 \times \left( t_0-\frac{\ep}{\sqrt{N}},
t_0 + \frac{\ep}{\sqrt{N}} \right)$.  Then, by (3.4) we have
$x \in \oooo\Omega$ and $d(x) > \frac{4}{N}
\Vert d\Vert_{C(\oooo{\Omega_1})}$, so that
$$
d(x) - \beta(t-t_0)^2
> \frac{4}{N}\Vert d\Vert_{C(\oooo{\Omega_1})}
- \frac{\beta\ep^2}{N},
$$
that is, $\va(x,t) > \mu_4$, which implies that $(x,t) \in D$
by the definition of $D$.
Next let $(x,t) \in D$.  Then $d(x) - \beta (t-t_0)^2
> \frac{1}{N}\Vert d\Vert_{C(\oooo{\Omega_1})}
- \frac{\beta \ep^2}{N}$.  Therefore
$$
\Vert d\Vert_{C(\oooo{\Omega_1})} - \frac{1}{N}
\Vert d\Vert_{C(\oooo{\Omega_1})} + \frac{\beta\ep^2}{N}
> \beta (t-t_0)^2.
$$
Applying (3.5), we have $2\left(1-\frac{1}{N}\right)\beta\ep^2
+ \frac{\beta\ep^2}{N} > \left(1 - \frac{1}{N}\right)\Vert d\Vert
_{C(\oooo{\Omega_1})} + \frac{\beta \ep^2}{N} > \beta(t-t_0)^2$, that is,
$2\beta \ep^2 > \beta(t-t_0)^2$, which implies that
$t_0 - \sqrt{2}\ep < t < t_0+\sqrt{2}\ep$.
The verification of (3.6) is completed.

Next we have
$$
\left\{
\begin{array}{rl}
&\ppp D \subset \Sigma_1 \cup \Sigma_2, \\
&\Sigma_1 \subset \Gamma \times (0,T), \quad
\Sigma_2 = \{ (x,t); \thinspace x\in\Omega, \thinspace
\va(x,t) = \mu_1\}.
\end{array}
\right.                           \eqno{(3.7)}
$$
In fact, let $(x,t) \in \ppp D$.  Then $x \in \oooo\Omega$ and
$\va(x,t) \ge \mu_1$.  We separately consider the cases
$x \in \Omega$ and $x \in \ppp\Omega$.  First let $x \in \Omega$.
If $\va(x,t) > \mu_1$, then $(x,t)$ is an interior point of
$D$, which is impossible.
Therefore $\va(x,t) = \mu_1$, which implies $(x,t) \in
\Sigma_2$.
Next let $x \in \ppp\Omega$.  Let $x \in \ppp\Omega \setminus
\Gamma$.  Then $x \in \ppp\Omega_1$ by the third condition in
(3.2), and $d(x) = 0$ by the second condition in (3.3).
On the other hand, $\va(x,t) \ge \mu_1$ yields that
$$
d(x) - \beta(t-t_0)^2 = -\beta(t-t_0)^2
\ge \frac{1}{N}\Vert d\Vert_{C(\oooo{\Omega_1})}
- \frac{\beta\ep^2}{N},
$$
that is, $0\le \beta(t-t_0)^2 \le \frac{1}{N}(
-\Vert d\Vert_{C(\oooo{\Omega_1})} + \beta\ep^2)$, which is impossible
by (3.5).  Therefore $x \in \Gamma$.  By (3.6), we see that
$0 < t < T$ and the verification of (3.7) is completed.

We apply Theorem 1 in $D$.
Henceforth $C>0$ denotes generic constants independent of $s$ and
choices of $v, p$.
We need a cut-off function because we have no data on
$\ppp D \setminus (\Gamma \times (0,T))$.
Let $\chi \in C^{\infty}(\R^{n+1})$ satisfying $0 \le \chi \le 1$ and
$$
\chi(x,t) =
\left\{
\begin{array}{rl}
1, \qquad & \va(x,t) > \mu_3, \\
0, \qquad & \va(x,t) < \mu_2.
\end{array}
\right.                             \eqno{(3.8)}
$$
We set $y = \chi v$ and $q = \chi p$.  Then, by (1.1) and (1.2), we have
\begin{align*}
& \ppp_ty - \kappa\Delta y + (A\cdot \nabla)y + (y\cdot \nabla)B
+ \nabla q\\
=& \chi F + v\ppp_t\chi - 2\kappa\nabla\chi\cdot \nabla v
- \kappa(\Delta\chi)v + (A\cdot\nabla \chi)v + p(\nabla\chi)
\quad \mbox{in $D$}
\end{align*}
and
$$
\ddd y = \nabla\chi \cdot v \quad \mbox{in $D$}.
$$
By (3.7) and (3.8), we see that
$$
\vert y\vert  = \vert \nabla y\vert = \vert q\vert = 0
\qquad \text{on $\Sigma_2$}.
$$
Hence Theorem 1 yields
\begin{align*}
&  \Vert (y,q)\Vert^2_{ \mathcal X_s(D)}
\le C\int_D \vert F\vert^2 \weight dxdt \\
+ &C\int_D  \vert v\ppp_t\chi - 2\kappa\nabla\chi\cdot \nabla v
- \kappa(\Delta\chi)v + (A\cdot\nabla \chi)v + p(\nabla\chi)\vert^2
\weight dxdt\\
+ &C\int_D (\vert \nabla\chi\cdot v\vert^2
+ \vert \nabla_{x,t}(\nabla\chi\cdot v)\vert^2)\weight dxdt
\end{align*}
$$
+ Ce^{Cs}(\Vert \chi v\Vert^2_{L^2(0,T;H^{\frac{3}{2}}(\Gamma))}
+ \Vert \ppp_t(\chi v)\Vert^2_{L^2(0,T;H^{\frac{1}{2}}(\Gamma))}
+ \Vert \sigma(\chi v, \chi p)\nu\Vert^2_{L^2(0,T;H^{\frac{1}{2}}(\Gamma))})
                                                     \eqno{(3.9)}
$$
for $s \ge s_0$.
We can verify $\Vert \chi v\Vert_{H^{\gamma}(\Gamma)}
\le C\Vert v\Vert_{H^{\gamma}(\Gamma)}$ with $\gamma = 0,1,2$, and for
$j=\frac{1}{2}$ and $j=\frac{3}{2}$, the interpolation inequality yields
$$
\Vert \chi v\Vert^2_{L^2(0,T;H^j(\Gamma))}
\le C\Vert v\Vert^2_{L^2(0,T;H^j(\Gamma))}, \quad
\Vert \ppp_t(\chi v)\Vert^2_{L^2(0,T;H^{\frac{1}{2}}(\Gamma))}
\le C\Vert \ppp_tv\Vert^2_{L^2(0,T;H^{\frac{1}{2}}(\Gamma))}.
$$
Therefore, since 
$$
\sigma(\chi v, \chi p)\nu = \chi\sigma(v,p)\nu
+ \kappa((\ppp_i\chi)v_j + (\ppp_j\chi)v_i)_{1\le i,j\le n}\nu,
$$
we have 
$$
\Vert \sigma(\chi v, \chi p)\nu\Vert
_{L^2(0,T;H^{\frac{1}{2}}(\Gamma))}
\le \Vert \sigma(v,p)\nu\Vert_{L^2(0,T;H^{\frac{1}{2}}(\Gamma))}
+ C\Vert v\Vert_{L^2(0,T;H^{\frac{1}{2}}(\Gamma))}
$$
by $\chi \in C^{\infty}(\R^{n+1})$.
Hence
\begin{align*}
&\Vert \chi v\Vert^2_{L^2(0,T;H^{\frac{3}{2}}(\Gamma))}
+ \Vert \ppp_t(\chi v)\Vert^2_{L^2(0,T;H^{\frac{1}{2}}(\Gamma))}
+ \Vert \sigma(\chi v, \chi p)\nu\Vert^2_{L^2(0,T;H^{\frac{1}{2}}(\Gamma))}\\
\le& C(\Vert v\Vert^2_{L^2(0,T;H^{\frac{3}{2}}(\Gamma))}
+ \Vert \ppp_tv\Vert^2_{L^2(0,T;H^{\frac{1}{2}}(\Gamma))}
+ \Vert \sigma(v, p)\nu\Vert^2_{L^2(0,T;H^{\frac{1}{2}}(\Gamma))}).
\end{align*}
We recall that
$$
G^2 = \Vert F\Vert^2_{L^2(Q)}
+ \Vert v\Vert^2_{L^2(0,T;H^{\frac{3}{2}}(\Gamma))}
+ \Vert \ppp_t v\Vert^2_{L^2(0,T;H^{\frac{1}{2}}(\Gamma))}
+ \Vert \sigma(v, p)\nu\Vert^2_{L^2(0,T;H^{\frac{1}{2}}(\Gamma))}.
$$

The integrands of the second and the third terms on the right-hand
side of (3.9) do not vanish only if $\va(x,t) \le \mu_3$, because
these coefficients include derivatives of $\chi$ as factors and
by (3.8) vanish if $\va(x,t) > \mu_3$.  Therefore
\begin{align*}
&\vert \mbox{[the second and the third terms on the right-hand side of
(3.9)]}
\vert\\
\le &C(\Vert v\Vert^2_{H^{1,1}(Q)}
+ \Vert p\Vert_{L^2(Q)}^2)e^{2s\mu_3}.
\end{align*}
Consequently (3.9) yields
$$
 \Vert (y,q)\Vert^2_{ \mathcal X_s(D)}
\le C(\Vert v\Vert^2_{H^{1,1}(Q)} + \Vert p\Vert_{L^2(Q)}^2)e^{2s\mu_3}
+ Ce^{Cs}G^2     \quad \forall s\ge s_0.                      \eqno{(3.10)}
$$

By (3.4) and the definition of $D$, we can directly verify that
$(x,t) \in \Omega_0 \times \left(t_0 - \frac{\ep}{\sqrt{N}},
t_0 + \frac{\ep}{\sqrt{N}}\right)$ implies $\va(x,t) > \mu_4$.
Therefore, noting (3.6) and (3.8), we see that
\begin{align*}
&  \Vert (y,q)\Vert^2_{ \mathcal X_s(D)}\ge  \Vert (v,p)\Vert^2_{ \mathcal X_s(\Omega_0\times (t_0 - \frac{\ep}{\sqrt{N}},t_0 + \frac{\ep}{\sqrt{N}}) )}\\
\ge& e^{2s\mu_4}\int^{t_0 + \frac{\ep}{\sqrt{N}}}_{t_0 - \frac{\ep}{\sqrt{N}}}
\int_{\Omega_0}
\left\{ \frac{1}{s^2}\left( \vert \ppp_tv\vert^2
+ \ssss \vert \ppp_i\ppp_jv \vert^2\right)
+ \vert \nabla v\vert^2 + s^2\vert v\vert^2
+ \frac{1}{s}\vert \nabla p\vert^2 + s\vert p\vert^2 \right\} dxdt.
\end{align*}
Hence (3.10) yields
\begin{align*}
& e^{2s\mu_4}\int^{t_0+\frac{\ep}{\sqrt{N}}}_{t_0-\frac{\ep}{\sqrt{N}}}
\int_{\Omega_0}
\left\{ \frac{1}{s^2}\left( \vert \ppp_tv\vert^2
+ \ssss \vert \ppp_i\ppp_jv \vert^2\right)
+ \vert \nabla v\vert^2 + s^2\vert v\vert^2
+ \frac{1}{s}\vert \nabla p\vert^2 + s\vert p\vert^2 \right\} dxdt\\
\le & C(\Vert v\Vert^2_{H^{1,1}(Q)} + \Vert p\Vert_{L^2(Q)}^2)e^{2s\mu_3}
+ Ce^{Cs}G^2.
\end{align*}
Therefore
\begin{align*}
& \int^{t_0+\frac{\ep}{\sqrt{N}}}_{t_0-\frac{\ep}{\sqrt{N}}}
\int_{\Omega_0}
\left\{ \left( \vert \ppp_tv\vert^2
+ \ssss \vert \ppp_i\ppp_jv \vert^2\right)
+ \vert \nabla v\vert^2 + \vert v\vert^2
+ \vert \nabla p\vert^2 + \vert p\vert^2 \right\}dxdt\\
\le & Cs^2e^{-2s(\mu_4-\mu_3)}
(\Vert v\Vert^2_{H^{1,1}(Q)} + \Vert p\Vert_{L^2(Q)}^2)
+ Ce^{Cs}G^2\quad\forall s\ge s_0.
\end{align*}

By $\sup_{s>0} se^{-s(\mu_4-\mu_3)} < \infty$, we
estimate $se^{-2s(\mu_4-\mu_3)}$ by $e^{-s(\mu_4-\mu_3)}$ on
the right-hand side.
Moreover, replacing $C$ by $Ce^{Cs_0}$, we can have
$$
\Vert v\Vert^2_{H^{2,1}\left(\Omega_0\times
\left(t_0-\frac{\ep}{\sqrt{N}},t_0+\frac{\ep}{\sqrt{N}}\right)\right)}
+ \Vert p\Vert^2_{H^{1,0}\left(\Omega_0\times
\left(t_0-\frac{\ep}{\sqrt{N}},t_0+\frac{\ep}{\sqrt{N}}\right)\right)}
$$
$$
\le Ce^{-s(\mu_4-\mu_3)}(\Vert v\Vert^2_{H^{1,1}(Q)}
+ \Vert p\Vert^2_{L^2(Q)}) + Ce^{Cs}G^2   \eqno{(3.11)}
$$
for all $s \ge 0$.
Let $m \in \N$ satisfy 
$\sqrt{2}\ep + \frac{m\ep}{\sqrt{N}} \le T - \sqrt{2}\ep
\le \sqrt{2}\ep + \frac{(m+1)\ep}{\sqrt{N}} \le T$.

We here notice that the constant $C$ in (3.11) is independent
also of $t_0$, provided that $\sqrt{2}\ep \le t_0 \le T-\sqrt{2}\ep$.
In (3.11), taking $t_0 = \sqrt{2}\ep + \frac{j\ep}{\sqrt{N}}$,
$j=0,1,2,..., m$ and summing up over $j$, we have
$$
\Vert v\Vert^2_{H^{2,1}\left(\Omega_0\times
\left(\sqrt{2}\ep- \frac{\ep}{\sqrt{N}}, T-\sqrt{2}\ep - \frac{\ep}{\sqrt{N}}
\right)\right)}
+ \Vert p\Vert^2_{H^{1,0}\left(\Omega_0\times
\left(\sqrt{2}\ep- \frac{\ep}{\sqrt{N}}, T-\sqrt{2}\ep - \frac{\ep}{\sqrt{N}}
\right)\right)}
$$
$$
\le Ce^{-s(\mu_4-\mu_3)}(\Vert v\Vert^2_{H^{1,1}(Q)}
+ \Vert p\Vert^2_{L^2(Q)})
+ Ce^{Cs}G^2
$$
for all $s \ge 0$.
Here we note that $T - \sqrt{2}\ep \le \sqrt{2}\ep + \frac{(m+1)\ep}
{\sqrt{N}}$ implies
$T - \sqrt{2}\ep - \frac{m\ep}{\sqrt{N}} \le \sqrt{2}\ep + \frac{1}{\sqrt{N}}
\ep$.  Replacing $\left(\sqrt{2}+\frac{1}{\sqrt{N}}\right)\ep$ by
$\ep$, we have
$$
\Vert v\Vert^2_{H^{2,1}(\Omega_0\times (\ep, T-\ep))}
+ \Vert p\Vert^2_{H^{1,0}(\Omega_0\times (\ep,T-\ep))}
$$
$$
\le Ce^{-s(\mu_4-\mu_3)}(\Vert v\Vert^2_{H^{1,1}(Q)}
+ \Vert p\Vert^2_{L^2(Q)})
+ Ce^{Cs}G^2                               \eqno{(3.12)}
$$
for all $s \ge s_0$.

First let $G=0$.  Then letting $s \to \infty$ in (3.12), we see that
$\vert v\vert = \vert p\vert =0$ in $\Omega_0 \times (\ep,T-\ep)$, so
that the conclusion of Theorem 2 holds true.  Next let $G\ne 0$.
First let $G \ge \Vert v\Vert_{H^{1,1}(Q)} +
\Vert p\Vert_{L^2(Q)}$.  Then (3.12) implies
$\Vert v\Vert_{H^{2,1}(\Omega_0\times (\ep,T-\ep))}
+ \Vert p\Vert_{H^{1,0}(\Omega_0\times (\ep,T-\ep))}
\le Ce^{Cs}G$ for $s \ge 0$, which already proves the theorem.
Second let $G < \Vert v\Vert_{H^{1,1}(Q)} + \Vert p\Vert
_{L^2(Q)}$.  In order to make the right-hand side of (3.12) smaller,
we choose $s > 0$ such that
$$
e^{-s(\mu_4-\mu_3)}(\Vert v\Vert^2_{H^{1,1}(Q)}
+ \Vert p\Vert^2_{L^2(Q)}) = e^{Cs}G^2.
$$
By $G\ne 0$, we can choose
$$
s = \frac{1}{C+\mu_4-\mu_3} \log\frac{\Vert v\Vert^2_{H^{1,1}(Q)}
+ \Vert p\Vert^2_{L^2(Q)}}{G^2} > 0.
$$
Then (3.12) gives
$$
\Vert v\Vert^2_{H^{2,1}(\Omega_0\times (\ep,T-\ep))}
+ \Vert p\Vert^2_{H^{1,0}(\Omega_0\times (\ep,T-\ep))}
\le 2C(\Vert v\Vert^2_{H^{1,1}(Q)}
+ \Vert p\Vert^2_{L^2(Q)})^{\frac{C}{C+\mu_4-\mu_3}}
G^{\frac{2(\mu_4-\mu_3)}{C+\mu_4-\mu_3}}.
$$
The the proof of Theorem 2 is completed. $\blacksquare$
\\

\end{document}